\documentclass[onecolumn]{autart}
\usepackage{amsmath}
\usepackage{amssymb}
\usepackage{color}

\usepackage{pgfplots}
\usetikzlibrary{intersections}
\usetikzlibrary{patterns}
\usetikzlibrary{positioning}
\usetikzlibrary{calc}
\usepgfplotslibrary{patchplots}

\newcommand{\tr}{\intercal}
\newcommand{\dist}[1]{\mathrm{dist}\left(#1\right)}

\begin{document}

\begin{frontmatter}

\title{
Towards Global Optimal Control via Koopman Lifts\thanksref{footnoteinfo}} 

\thanks[footnoteinfo]{MEV and BH acknowledge financial support via ShanghaiTech U.~Grant F-0203-14-012. CJ would like to acknowledge support received from the Swiss National Science Foundation under the RISK project (Risk Aware Data-Driven Demand Response), grant number 200021 175627.}

\author[Shanghai]{Mario E.~Villanueva}\ead{meduardov@shanghaitech.edu.cn},
\author[EPFL]{Colin N.~Jones}\ead{colin.jones@epfl.ch},
\author[Shanghai]{Boris Houska}\ead{borish@shanghaitech.edu.cn}

\address[Shanghai]{School of Information Science and Technology, ShanghaiTech, China}%

\address[EPFL]{Automatic Control Laboratory, EPFL, Switzerland}%


\begin{abstract}
This paper introduces a framework for solving time-autonomous nonlinear infinite horizon optimal control problems, under the assumption that all minimizers satisfy Pontryagin's necessary optimality conditions. In detail, we use methods from the field of symplectic geometry to analyze the eigenvalues of a Koopman operator that lifts Pontryagin's differential equation into a suitably defined infinite dimensional symplectic space. This has the advantage that methods from the field of spectral analysis can be used to characterize globally optimal control laws. A numerical method for constructing optimal feedback laws for nonlinear systems is then obtained by computing the eigenvalues and eigenvectors of a matrix that is obtained by projecting the Pontryagin-Koopman operator onto a finite dimensional space. We illustrate the effectiveness of this approach by computing accurate approximations of the optimal nonlinear feedback law for a Van der Pol control system, which cannot be stabilized by a linear control law.
\end{abstract}

\end{frontmatter}

\section{Introduction}
During the last decades algorithms and software for computing local solutions of nonlinear optimal control problems have reached a high level of maturity~\cite{Biegler2007,Diehl2002,Zavala2009}, enabling its deployment in industrial applications~\cite{Qin2003}. By now, such local solutions can be computed within the milli- and microsecond range. Moreover, auto-generated implementations of real-time local optimal control solvers run on embedded hardware systems~\cite{Mattingley2009,Houska2011}. In contrast to these developments in local optimal control, algorithms for locating global minimizers of non-convex optimal control problems, can hardly be claimed to be ready for widespread industrial application. There are at least two reasons for this. Firstly, generic global optimization methods can often only be applied to problems of modest dimensions. And secondly, their run-time usually exceeds the run-time of local solvers by orders of magnitude. Nevertheless, in the last years there have been a number of promising developments in the field of global optimal control, which are reviewed next.

Dynamic Programming (DP)~\cite{Bellman1957}, which proceeds by approximately solving the Hamilton-Jacobi-Bellman (HJB) equation~\cite{Frankowska1993}, is a historically important method able to find globally optimal feedback laws. In~\cite{Luus1990} and~\cite{Gruene2002} tailored discretization grids for DP implementations have been developed that can successfully solve nonlinear optimal control problems, as long as the dimension of the state of the system is small. For higher dimensional state-spaces these methods are, however, not applicable as the complexity of storing and processing the value functions during the DP recursion grows exponentially with the number of states. Other deterministic methods for global optimization problems involving nonlinear differential equations are based on Branch-and-Bound (BB)~\cite{Chachuat2006,Esposito2000,Lin2007b,Singer2006} and its variants, including \mbox{$\alpha$-BB} methods~\cite{Papamichail2002,Diedam2017}. These BB algorithms have in common that they can effectively solve problems with a small number of decision variables. Unfortunately, implementations of branch-and-bound search easily run out of memory in higher dimensional spaces due to exponentially growing search trees. In particular, as the discretization of control inputs typically leads to a large number of optimization variables, BB-methods are usually not suited for solving optimal control problems.

An alternative to Branch-and-Bound are the so-called Branch-and-Lift (BL) methods~\cite{Houska2014}. In contrast to BB, these methods never discretize the control inputs directly.
Instead, these methods branch over orthogonal projections of the control function in lower dimensional subspaces of increasing dimension, until an $\epsilon$-suboptimal global control input is found. A rigorous analysis of the mathematical properties of such BL methods can be found in~\cite{Houska2019}.
However, in the context of developing practical implementations of Branch-and-Lift one faces numerical challenges. In particular, these methods require the computation of accurate enclosures of moment-constrained reachable sets of nonlinear differential equations. At the current state of research, the lack of generic methods for computing such moment-constrained reachable set enclosures limits the applicability of Branch-and-Lift---currently, some problems have been solved successfully with BL by using enclosure propagation methods that are tailored to particular applications~\cite{Feng2017,Houska2014}.

In order to mitigate the above-mentioned limitations of existing global optimal control
methods, this paper proposes a completely new framework for constructing methods for solving nonlinear infinite horizon optimal control problems. Here, the main idea is to analyze the spectrum of a Koopman operator that is associated with Pontryagin's differential equation---under the assumption that Pontryagin's optimality condition~\cite{Pontryagin1962} is satisfied at the optimal solution. In order to understand the contributions of this paper, it is important to first be aware of the fact that Koopman operators
have originally been introduced almost $90$ years ago---in fact, by Koopman himself~\cite{Koopman1931}. This theory has later been extended for more general nonlinear differential equations in~\cite{Mezic2005}, where the concept of Koopman mode analysis has been introduced. A more complete analysis of these Koopman mode-based methods from the field of ergodic theory has, however, only appeared much later~\cite{Arbabi2017}. Notice that a mathematical analysis of finite dimensional approximations of the Koopman operator can be found in~\cite{Govindarajan2019} and a variety of related applications of approximate Koopman mode analysis methods can be found in~\cite{Budisic2012}. Moreover, in~\cite{Korda2018} Koopman mode estimation heuristics for data driven control have been introduced.

\paragraph*{Contribution.}
The key idea of this paper is to introduce a Koopman operator that is associated with the Pontryagin differential equation 
of a nonlinear infinite-horizon optimal control problem---in the following we called a 
\textit{Pontryagin-Koopman operator}---which can be used to construct globally optimal feedback laws. Section~\ref{sec::ocp} introduces infinite horizon optimal control problems and briefly reviews Pontryagin's necessary condition for optimality. This review section is followed by a presentation the main contributions of this article, which can be outlined as follows.

\begin{enumerate}

\item Section~\ref{sec::3} uses ideas from the field of symplectic geometry~\cite{Arnold2001} to characterize the structural properties of Pontryagin-Koopman operators. Section~\ref{sec::spectrum} leverages on these structural properties in order to perform a symplectic spectral analysis that finally leads to a complete characterization of globally optimal feedback control laws (see Section~\ref{sec::Manifolds}, Theorem~\ref{thm::OptimalControl}). Because this is the first time that the symplecticity of Pontryagin-Koopman operators is analyzed and used to derive practical characterizations of optimal control laws, Theorem~\ref{thm::OptimalControl} can be regarded as the main theoretical contribution of this paper.

\item In order to avoid misunderstandings, it is mentioned clearly that this paper does not claim to solve all the numerical issues regarding the discretization of Pontryagin-Koopman operators that would have to be solved to develop a generic software for global optimal control. However, Section~\ref{sec::numerics} illustrates the practical applicability of the proposed theoretical framework by designing an optimal regulator for a controlled Van der Pol oscillator. Although this implementation is only based on a naive Galerkin projection of the Pontryagin-Koopman operator onto a finite dimensional Legendre basis, this section successfully applies the proposed framework to construct accurate approximations of nonlinear globally optimal feedback \mbox{laws---in} this example, for a system that cannot be stabilized by a linear control law.

\end{enumerate}

The above contributions are relevant for the future of optimal control algorithm development, as they pave the way towards the development of practical procedures for approximating globally optimal control laws of a very general class of nonlinear systems. Therefore, Section~\ref{sec::conclusion} does not only conclude the paper, but also outlines the potential and relevance of the proposed framework for future research in control.

\paragraph*{Notation.} The distance of a point $x \in \mathbb R^{n_x}$ to a trajectory $\phi: \mathbb R \to \mathbb R^n$ is denoted by
\[
\dist{x,\phi} = \inf_{t \in \mathbb R} \Vert x - \phi(t) \Vert \; .
\]
We denote with $\mathbb L_p^n$ the set of (potentially complex-valued) Lebesgue measurable functions $\varphi: \mathbb R^{n} \to \mathbb C$ whose $p$-th power of the absolute value is integrable on $\mathbb R^{n}$. Moreover, we use the notation $\mathbb W_{k,p}^n$ to denote the Sobolev space of integrable functions on $\mathbb R^n$, whose weak derivatives up to order $k$ are all in $\mathbb L_p^n$. The symbol $M^\tr$ denotes the Hermitian transpose of the matrix 
$M \in \mathbb C^{n \times m}$. The symbol $\nabla$ denotes the gradient operator. The associated second order derivative operator is denoted by $\nabla^2 = \nabla \nabla^\tr$. Last but not least, the support of a function $\varphi: \mathbb R^n \to \mathbb C$ is denoted by
$$\mathrm{supp}(\varphi) = \{ x \in \mathbb R^{n} \mid \varphi(x) \neq 0 \} \; .$$

\section{Infinite horizon optimal control}
\label{sec::ocp}
This section introduces infinite horizon optimal control problems and briefly summarizes existing methods for analyzing the local stability properties of optimal periodic limit orbits (see Section~\ref{subsec::periodic}). Moreover, Sections~\ref{subsec::Pontryagin} and~\ref{subsec::Boundary} review Pontryagin's necessary optimal condition thereby introducing the notation that is used throughout the paper.

\subsection{Problem formulation}
This paper considers infinite horizon optimal control problems of the form
\begin{equation}
\label{eq::ocp}
\begin{alignedat}{2}
V(x_0) = &\min_{x,u} && \displaystyle\int_{0}^{\infty} \left[ l(x(t),u(t)) - l^\star \right] \, \mathrm{d}t \\[0.3cm]
&\mathrm{s.t.}&& \left\{
\begin{aligned}
&\forall t \in [0,\infty),\\[0.1cm]
&\dot x(t) = f(x(t),u(t)) \\[0.1cm]
&x(0) = x_0 \; .
\end{aligned}
\right.
\end{alignedat}
\end{equation}
Here, $x: \mathbb R \to \mathbb R^{n_x}$ denotes the state trajectory and 
$u: \mathbb R \to \mathbb R^{n_u}$ denotes the control input.
The initial value $x_0 \in \mathbb R^{n_x}$ is assumed to be given. Throughout this paper the following assumptions are imposed.

\begin{assum}
\label{ass::f}
The function $f: \mathbb R^{n_x} \times \mathbb R^{n_u} \to \mathbb R^{n_x}$ is twice continuously differentiable and globally Lipschitz continuous in $x$.
\end{assum}

\begin{assum}
\label{ass::l}
The function $l: \mathbb R^{n_x} \times \mathbb R^{n_u} \to \mathbb R$ is twice continuously differentiable.
\end{assum}
The constant $l^\star \in \mathbb R$ in~\eqref{eq::ocp} denotes the optimal average cost,
\begin{align}
l^\star = \lim_{T \to \infty} & \min_{x,u} \frac{1}{T} \int_{0}^T l(x(t),u(t)) \, \mathrm{d}t \notag \\
& \; \, \mathrm{s.t.} \; \; \left\{
\begin{aligned}
&\forall t \in [0,T],\\[0.1cm]
&\dot x(t) = f(x(t),u(t)) \\[0.1cm]
&x(0) = x_0 \; ,
\end{aligned}
\right. \notag
\end{align}
assuming that this limit exists. Notice that $l^\star$ does not need to be known explicitly in order to solve~\eqref{eq::ocp}, as adding constant offsets to the objective does not affect the solution of an optimization problem. In this paper, the constant $l^\star$ is merely introduced for mathematical reasons, such that $V(x_0)$ remains finite and well-defined on infinite horizons under mild regularity assumptions that shall be introduced further below.

\subsection{Limit behavior and periodic orbits}
\label{subsec::periodic}

In practical instances, one is often interested in whether optimal solutions for the state trajectory of~\eqref{eq::ocp} converge to an optimal steady-state or, in more generality,
to an optimal periodic limit orbit. These optimal steady states or more general periodic orbits are defined as follows.

\begin{defn}
The function $(x_\mathrm{p},u_\mathrm{p}) \in \mathbb W_{1,1}^{n_x} \times \mathbb L_1^{n_x}$ is called an optimal periodic orbit, if there exists a period $T > 0$ such that for all $t \in [0,\infty)$
\begin{enumerate}
\addtolength{\itemsep}{4pt}
\item $x_\mathrm{p}(t+T) = x_\mathrm{p}(t)$ and $u_\mathrm{p}(t+T) = u_\mathrm{p}(t)$,
\item $\dot x_\mathrm{p}(t) = f(x_\mathrm{p}(t),u_\mathrm{p}(t))$ for all $t \in [0,T]$, and
\item $\frac{1}{T} \int_0^{T} l(x_\mathrm{p}(t),u_\mathrm{p}(t)) \, \mathrm{d}t = l^\star$.
\end{enumerate}
For the special case that $x_\mathrm{p}$ and $u_\mathrm{p}$ are constant functions, this pair is called an optimal steady-state.
\end{defn}

In order to prepare the following analysis, it is helpful to introduce the shorthands
\begin{equation}
\begin{aligned}
A(t) &= f_x(x_\mathrm{p}(t),u_\mathrm{p}(t)) \;, \quad
B(t)  = f_u(x_\mathrm{p}(t),u_\mathrm{p}(t)) \;, \\[0.1cm]
Q(t) &= l_{xx}(x_\mathrm{p}(t),u_\mathrm{p}(t)) \;, \quad 
R(t)  = l_{uu}(x_\mathrm{p}(t),u_\mathrm{p}(t)) \;,  \\[0.1cm]
S(t) &= l_{xu}(x_\mathrm{p}(t),u_\mathrm{p}(t))
\end{aligned} 
\end{equation}
assuming that $(x_\mathrm{p},u_\mathrm{p})$ is an optimal periodic orbit. Here, $f_x$ and $f_u$ denote the partial derivatives of $f$ with respect to $x$ and $u$ and an analogous notation is then also used for the mixed second order derivatives of $l$. The above matrix-valued
functions can be used to construct sufficient conditions under which an optimal
periodic orbit is locally stabilizable. Here, one relies on the theory of periodic Riccati differential equations~\cite{Bittanti1991} of the form
\begin{equation}
\begin{aligned}
-\dot P(t) &= P(t)A(t) + A(t)^\tr P(t) + Q(t) 
- (P(t) B(t) + S(t))R(t)^{-1} (P(t) B(t) + S(t))^\tr \\[0.1cm]
\label{eq::Riccati}
P(0) &= P(T) \succ 0 \; .
\end{aligned}
\end{equation}
It is well-known~\cite{Bittanti1991} that if
Assumptions~\ref{ass::f} and~\ref{ass::l} are satisfied and if $R$ has full rank, then the existence of a periodic and positive definite function $P$ satisfying~\eqref{eq::Riccati} is sufficient to ensure that the solution $x^\star$ of~\eqref{eq::ocp} converges to the optimal periodic orbit $x_\mathrm{p}$ as long as $\, \dist{x_0,x_\mathrm{p}} \,$ is sufficiently small---that is, if $x_0$ is in a small neighborhood of the optimal orbit $x_\mathrm{p}$. However, this statement is of a rather local nature. This means that, if we wish to understand the global behavior of system~\eqref{eq::ocp}, an analysis of the periodic Riccati equation is not sufficient. The 
following sections focus on analyzing global solutions of~\eqref{eq::ocp}, under the assumption that $\, \dist{x_0,x_\mathrm{p}} \,$ is not necessarily small.

\subsection{Pontryagin's differential equations}
\label{subsec::Pontryagin}

Pontryagin's maximum principle~\cite{Pontryagin1962} can be used to derive 
necessary conditions for the minimizers of~\eqref{eq::ocp}. The first 
order variational optimality condition is summarized as follows. 
Let $H: \mathbb R^{n_x} \times \mathbb R^{n_u} \times \mathbb R^{n_x} \to \mathbb R$
be the Hamiltonian function of~\eqref{eq::ocp},
\begin{align}
\label{eq::H}
H(x,u,\lambda) = \lambda^\tr f(x,u) + l(x,u) \; .
\end{align}
If Assumptions~\ref{ass::f} and~\ref{ass::l} hold, $H$ is, by construction, twice continuously differentiable in all variables. Next, let
\begin{align}
\label{eq::u}
u^\star(x,\lambda) = \underset{u}{\text{argmin}} \; H(x,u,\lambda)
\end{align}
denote the associated parametric minimizer of $H$. At this point, we introduce the following regularity assumption.
\begin{assum}
\label{ass::u}
The parametric minimizer $u^\star$ in~\eqref{eq::u} satisfies the second order sufficient condition, $H_{uu} \succ 0$.
\end{assum}
Notice that for the practically relevant special case that $f$ is affine in $u$ and $l$ strongly convex in $u$, Assumption~\eqref{ass::u} always holds whenever Assumptions~\ref{ass::f} and~\ref{ass::l} are satisfied.

Next, if Assumptions~\ref{ass::f},~\ref{ass::l}, and~\ref{ass::u} are satisfied, 
it is well-known~\cite{Pontryagin1962} that any minimizer $(x,u)$ of~\eqref{eq::ocp} necessarily satisfies Pontryagin's differential equation
\begin{align}
\label{eq::PontryaginODE1}
\dot x(t) &= f( x(t), u^\star(x(t),\lambda(t)) ) \\[0.16cm]
\label{eq::PontryaginODE2}
\dot \lambda(t) &= -\nabla_x H(x(t),u^\star(x(t),\lambda(t)),\lambda(t))
\end{align}
for a co-state function $\lambda \in \mathbb W_1^{n_x}$. In the following, we summarize this differential equation in the form
\begin{align}
\label{eq::PontryaginODE}
\dot y(t) = F(y(t)) \; ,
\end{align}
where $y = [ x^\tr, \lambda^\tr ]^\tr$ denotes the stacked state and $F$ a stacked version of the right-hand side functions of~\eqref{eq::PontryaginODE1} and~\eqref{eq::PontryaginODE2}.

\subsection{Necessary boundary and limit conditions}
\label{subsec::Boundary}
Besides Pontryagin's differential equation~\eqref{eq::PontryaginODE}, optimal solutions of~\eqref{eq::ocp} satisfy necessary boundary conditions. First, of  course, the initial condition $x(0) = x_0$ must hold. Moreover, the co-state $\lambda$ satisfies a necessary limit condition, which can be summarized as follows.

\begin{prop}
\label{prop::multiplier}
Let Assumptions~\ref{ass::f},~\ref{ass::l}, and~\ref{ass::u} hold.
Moreover, let $(x^\star,u^\star)$ be a primal optimal solution of~\eqref{eq::ocp}
converging to an optimal periodic orbit $(x_\mathrm{p},u_\mathrm{p}) \in \mathbb W_{1,1}^{n_x} \times \mathbb L_1^{n_x}$ at which the
periodic Riccati differential equation~\eqref{eq::Riccati} admits a positive definite solution. Then, there exists a periodic function $\lambda_{\mathrm{p}} \in \mathbb W_{1,1}^{n_x}$ such
that the associated co-state $\lambda^\star$ necessarily satisfies
\[
\lim_{t \to \infty} \Vert \lambda^\star(t) - \lambda_{\mathrm{p}}(t) \Vert = 0 \; .
\]
\end{prop}

Notice that Proposition~\ref{prop::multiplier} is---at least in very similar versions---well-known in the literature~\cite{Pontryagin1962,Liberzon2012}. However, as this result is important for understanding the developments in this paper, the following proof briefly recalls the main argument, why this proposition holds.

\textbf{Proof.}
Since~\eqref{eq::Riccati} admits a positive definite solution, the value function
$V$ in~\eqref{eq::ocp} is well-defined and differentiable in a neighborhood of the
optimal orbit $x^\star$. Thus,~\eqref{eq::ocp} can be replaced by an equivalent 
finite-horizon optimal control problem as long as $V$ is used as a terminal cost.
Now, Pontryagin's principle for optimal control problems with Mayer 
terms~\cite{Pontryagin1962} yields the boundary condition
\[
\lambda^\star(t) = \nabla V(x^\star(t))
\]
for any horizon length $t > 0$. Clearly, since $x^\star$ converges to the optimal periodic orbit $x_\mathrm{p}$, $\lambda^\star$ must converge to the associated dual limit orbit $
\lambda_{\mathrm{p}}(t) = \nabla V(x_{\mathrm{p}}(t))$.
\qed

In the next sections we will see that Proposition~\ref{prop::multiplier} implies that the states of Pontryagin's differential equation for stabilizable systems must evolve along certain stable manifolds. As it shall demonstrated, these manifolds can be characterized using 
a Koopman mode analysis~\cite{Mezic2005}.

\section{Symplectic Koopman operators}
\label{sec::3}
This section analyzes the symplectic structure of the flow associated with Pontryagin's differential equation. These symplectic structures are needed to understand the properties of the Pontryagin Koopman operators. In fact, Section~\ref{sec::Koopman} introduces a symplectic test space---that is, a space of observables---in which the Pontryagin-Koopman operator inherits certain symplecticity properties of its underlying incompressible flow field. Moreover, Section~\ref{sec::duality} leverages on ideas from the field of symplectic geometry~\cite{Arnold2001} in order to work out the symplectic dual of the Pontryagin-Koopman operator. Notice that these developments are the basis for the developments in Section~\ref{sec::spectrum}, which uses a symplectic spectral analysis in order to characterize globally optimal control laws.

\subsection{Symplectic Flows}
\label{sec::flow}
Let $\Gamma_t: \mathbb R^{2n_x} \to \mathbb R^{2n_x}$ denote the flow
of Pontryagin's differential equations such that
\[
\frac{\mathrm{d}}{\mathrm{d}t} \Gamma_t(z) = F( \Gamma_t(z) ) \quad \text{and} \quad \Gamma_0(z) = z
\]
for all $z \in \mathbb R^{2n_x}$.
If Assumptions~\ref{ass::f},~\ref{ass::l}, and~\ref{ass::u} hold, then
$\Gamma_t$ is a well-defined, continuously differentiable function. In the mathematical literature~\cite{Arnold2001} it is well-known that $\Gamma_t$ is a so-called symplectic flow. In order to reveal this symplectic structure it is helpful to introduce the block matrix
$$\Omega = \left(
\begin{array}{rr}
0 \; \, & I \\
-I \; \, & 0
\end{array}
\right) \in \mathbb R^{2n_x \times 2 n_x} \; .
$$
Now, the structural properties of $\Gamma_t$ can be summarized as follows.

\begin{lem}
\label{lem::symplectic}
Let Assumptions~\ref{ass::f},~\ref{ass::l}, and~\ref{ass::u} be satisfied. Then, the
function $\frac{\partial}{\partial z} \Gamma_t$ satisfies the equation
\begin{align}
\label{eq::Gsymplectic}
\forall t \in \mathbb R, \quad \left[ \frac{\partial}{\partial z} \Gamma_t \right]^\tr \Omega \left[ \frac{\partial}{\partial z} \Gamma_t \right] = \Omega \; ,
\end{align}
that is, $\Gamma_t$ is symplectic function.
\end{lem}

\textbf{Proof.} Let us first recall that the right-hand side $F$ of Pontryagin's differential equation is given by~\eqref{eq::PontryaginODE1} and~\eqref{eq::PontryaginODE2}, which
can also be summarized as
\begin{align}
F = \left[
\begin{array}{c}
f \\ -H_x^\tr
\end{array}
\right] \quad \text{with} \quad y = \left[
\begin{array}{c}
x \\ \lambda
\end{array}
\right] \; .
\end{align}
This implies that the derivative of $F$ with respect $y$ can be written in the form
\begin{align}
\label{eq::Fy}
F_{y} = \left(
\begin{array}{cc}
f_x - f_u H_{uu}^{-1} H_{ux} \; \; & - f_u H_{uu}^{-1} f_u^\tr \\
-H_{xx} + H_{xu} H_{uu}^{-1} H_{ux} \; \; & -f_x^\tr + H_{xu} H_{uu}^{-1} f_u^\tr 
\end{array}
\right) .
\end{align}
Here, we have used Assumptions~\ref{ass::f} and~\ref{ass::l} to ensure that the second derivatives of $H$ with respect to $x$ and $u$ exist. Moreover, we have used Assumption~\ref{ass::u}, which implies that $H_{uu}$ is invertible and, as a consequence, that the implicit function theorem holds. In particular, we have
$$u_x^\star = -H_{uu}^{-1} H_{ux} \quad \text{and} \quad u_{\lambda}^\star = -H_{uu}^{-1} H_{u \lambda} = -H_{uu}^{-1} f_u^\tr \; .$$
In this form, it becomes clear that the matrix $\Omega F_y$ is symmetric and we arrive at the intermediate result
\begin{align}
\label{eq::FHamiltonian}
\Omega F_y = F_y^\tr \Omega^\tr = - F_y^\tr \Omega \; .
\end{align}
In order to proceed, we write the first order variational differential equation for $\Gamma_t$ in the form 
\begin{align}
\label{eq::SecondVariation}
\frac{\mathrm{d}}{\mathrm{d}t} \left[ \frac{\partial}{\partial z} \Gamma_t \right] = F_y \left[ \frac{\partial}{\partial z} \Gamma_t \right] \quad \text{with} \quad \left[ \frac{\partial}{\partial z} \Gamma_0 \right] = I \; .
\end{align}
Now, the main idea of this proof is to show that the function
\begin{align}
\Delta(t) = \left[ \frac{\partial}{\partial z} \Gamma_t \right]^\tr \Omega \left[ \frac{\partial}{\partial z} \Gamma_t \right] - \Omega
\end{align}
vanishes, $\Delta(t) = 0$. Here, we first have $\Delta(0) = 0$ by construction, since $\left[ \frac{\partial}{\partial z} \Gamma_0 \right] = I$. Moreover, the derivative of $\Delta$ with respect to time is given by
\begin{align}
\dot \Delta(t) &= \frac{\mathrm{d}}{\mathrm{d}t} \left[ \frac{\partial}{\partial z} \Gamma_t \right]^\tr \Omega \left[ \frac{\partial}{\partial z} \Gamma_t \right] + \left[ \frac{\partial}{\partial z} \Gamma_t \right]^\tr \Omega \frac{\mathrm{d}}{\mathrm{d}t} \left[ \frac{\partial}{\partial z} \Gamma_t \right] \notag \\[0.16cm]
&\overset{\eqref{eq::SecondVariation}}{=} \left[ \frac{\partial}{\partial z} \Gamma_t \right]^\tr \left[ F_y^\tr \Omega + \Omega F_y \right] \left[ \frac{\partial}{\partial z} \Gamma_t \right] \overset{\eqref{eq::FHamiltonian}}{=} 0
\notag
\end{align}
Thus, in summary, we must have $\Delta(t) = 0$ for all $t \in \mathbb R$, which yields the statement of this lemma.
\qed

The following corollary summarizes two immediate consequences of Lemma~\ref{lem::symplectic} which are both equivalent to stating that $\Gamma_t$ is an incompressible flow.

\begin{cor}
\label{cor::det}
Let Assumptions~\ref{ass::f},~\ref{ass::l}, and~\ref{ass::u} be satisfied. Then, 
$\Gamma_t$ is an incompressible flow; that is, we have
\begin{align}
\label{eq::det}
\forall t \in \mathbb R, \qquad \mathrm{det}\left( \left[ \frac{\partial}{\partial z} \Gamma_t \right] \right) = 1 \; .
\end{align}
Moreover, the divergence of the associated vector field $F$ vanishes, $\mathrm{div}(F) = \nabla^\tr F = 0$.
\end{cor}

\textbf{Proof.}  By taking the determinant on both sides of~\eqref{eq::Gsymplectic} and using that $\mathrm{det}(\Omega) = 1$, we find that
\[
\mathrm{det}\left( \left[ \frac{\partial}{\partial z} \Gamma_t \right] \right)^2 = 1 \quad \Longleftrightarrow \quad \mathrm{det}\left( \left[ \frac{\partial}{\partial z} \Gamma_t \right] \right) = \pm 1 \; .
\]
Since $\frac{\partial}{\partial z} \Gamma_0 = I$ and since $\frac{\partial}{\partial z} \Gamma_t$ is continuous, this is only possible if~\eqref{eq::det} holds. Next, by taking the logarithm on both sides of~\eqref{eq::det} and differentiating with respect to time, we find\footnote{Alternatively, the equation $\mathrm{Tr}(F_y) = 0$ can also directly be found by substituting the explicit expression~\eqref{eq::Fy} from the proof of Lemma~\ref{lem::symplectic}.}
\begin{align}
0 &= \frac{\mathrm{d}}{\mathrm{d} t} \log \left( \mathrm{det}\left( \left[ \frac{\partial}{\partial z} \Gamma_t \right] \right) \right) \notag \\[0.16cm]
&= \mathrm{Tr} \left( \left[ \frac{\partial}{\partial z} \Gamma_t \right]^{-1} \frac{\mathrm{d}}{\mathrm{d}t} \left[ \frac{\partial}{\partial z} \Gamma_t \right] \right) \overset{\eqref{eq::SecondVariation}}{=} \mathrm{Tr}( F_y ) = \mathrm{div}(F) \; . \notag
\end{align}
This completes the proof of the corollary. \qed

\subsection{Pontryagin-Koopman operators}
\label{sec::Koopman}

The developments from the previous section can be used to analyze the structural properties of the Pontryagin-Koopman operator, which are defined to be the Koopman operator that is associated with the flow $\Gamma_t$ of Pontryagin's differential equation, or, more formally:

\begin{defn}[Pontryagin-Koopman Operator]
\label{def::PontryaginKoopman}
We use the notation $U_t: \mathbb W_{1,2}^{2n_x} \to \mathbb W_{1,2}^{2n_x}$ to denote the Pontryagin Koopman operator of $\Gamma_t$, which is defined such that
\[
\forall \Phi \in \mathbb W_{1,2}^{2n_x}, \ \forall t \in \mathbb R, \qquad U_t \Phi = \Phi \circ \Gamma_t
\]
with $\circ$ denoting the composition operator. 
\end{defn}

It is well known~\cite{Koopman1931,Mezic2005} that $U_t$ is a linear operator satisfying $U_{t_1+t_2} = U_{t_1}U_{t_2}$ for all
$t_1,t_2 \in \mathbb R$. Moreover, $U_t$ satisfies $U_{t}^{-1} = U_{-t}$, which follows by substituting $t=t_1=-t_2$ in the previous equation and using that $\Gamma_0 = \mathrm{id}$.

In order to introduce a notion of symplecticity in the space of observables of the Pontryagin-Koopman operator, we introduce the bilinear form
$$\omega: \mathbb W_{1,2}^{2n_x} \times \mathbb W_{1,2}^{2n_x} \to \mathbb R \; ,$$
which is defined as
\[
\omega( \varphi, \Phi ) = \int_{\mathbb R^{2n_x}} \nabla \varphi(z)^\tr \Omega \nabla \Phi(z) \, \mathrm{d}z 
\]
for all $\varphi,\Phi \in \mathbb W_{1,2}^{2n_x}$. Because we have
$\Omega^\tr = -\Omega$, the bilinear form $\omega$ is skew-symmetric,
\[
\omega( \varphi, \Phi ) = -\omega( \Phi, \varphi )  ,
\]
and $\left( \mathbb W_{1,2}^{2n_x}, \omega \right)$ is a symplectic space~\cite{Arnold2001}.

\begin{thm}
\label{thm::symplectic}
Let Assumptions~\ref{ass::f},~\ref{ass::l}, and~\ref{ass::u} be satisfied. Now, $U_t$ is for all $t \in \mathbb R$ a symplectic operator in the space $\left( \mathbb W_{1,2}^{2n_x}, \omega \right)$. This means that we have
\[
\omega( U_t \varphi, U_t \Phi ) = \omega( \varphi, \Phi )
\]
for all $\varphi,\Phi \in \mathbb W_{1,2}^{2n_x}$.
\end{thm}

\textbf{Proof.}
Assumptions~\ref{ass::f},~\ref{ass::l}, and~\ref{ass::u} imply that the function $\Gamma_t$ and its inverse $\Gamma_t^{-1} = \Gamma_{-t}$ are both continuously differentiable. This implies in particular that the equivalence
\[
\Phi \in  \mathbb W_{1,2}^{2n_x} \quad \Longleftrightarrow \quad \left[ \Phi \circ \Gamma_t \right] \in \mathbb W_{1,2}^{2n_x}
\]
holds. Next, the definition of the operator $U_t$ and the chain rule for differentiation imply that
\begin{align}
\label{eq::aux2}
\nabla [ U_t \Phi ] = \nabla [ \Phi \circ \Gamma_t ] =   \left[ \frac{\partial}{\partial z}\Gamma_t \right]^\tr (\nabla \Phi) \circ \Gamma_t
\end{align}
for all $\Phi \in  \mathbb W_{1,2}^{2n_x}$. Furthermore, because Assumptions~\ref{ass::f},~\ref{ass::l}, and~\ref{ass::u} hold, it follows from Lemma~\ref{lem::symplectic} that
\begin{align}
\left[ \frac{\partial}{\partial z} \Gamma_t \right]^\tr \Omega \left[ \frac{\partial}{\partial z} \Gamma_t \right] = \Omega \; . \notag
\end{align}
Let us multiply this equation with $\left[ \frac{\partial}{\partial z} \Gamma_t \right] \Omega$ from the left and substitute $\Omega^2 = -I$. This yields
\[
\left[ \frac{\partial}{\partial z} \Gamma_t \right] \Omega \left[ \frac{\partial}{\partial z} \Gamma_t \right]^\tr \Omega \left[ \frac{\partial}{\partial z} \Gamma_t \right] = - \left[ \frac{\partial}{\partial z} \Gamma_t \right] \; .
\]
Next, we multiply the latter equation with $\left[ \frac{\partial}{\partial z} \Gamma_t \right]^{-1} \Omega$ from the right and use the relation $\Omega^2 = -I$ once more, in order to arrive at the equation
\begin{align}
\label{eq::aux3}
\left[ \frac{\partial}{\partial z} \Gamma_t \right] \Omega \left[ \frac{\partial}{\partial z}\Gamma_t \right]^\tr = \Omega \; .
\end{align}
Thus, by substituting the previous equations, we get
\begin{equation}
\omega( U_t \varphi, U_t \Phi ) = \int_{\mathbb R^{2n_x}} \nabla [ U_t \varphi ](z)^\tr \Omega \nabla [U_t \Phi](z) \, \mathrm{d}z \notag 
\overset{\eqref{eq::aux2},\eqref{eq::aux3}}{=} \int_{\mathbb R^{2n_x}} \nabla \varphi( \Gamma_t(z) )^\tr \Omega \nabla \Phi( \Gamma_t(z) ) \, \mathrm{d}z \, . \notag 
\label{eq::aux4}
\end{equation}
In order to further transform this integral, we need to introduce the change of variables $z' = \Gamma_t(z)$. Because Corollary~\ref{cor::det} ensures that
\[
\left| \mathrm{det} \left( \frac{\partial}{\partial z} \Gamma_t(z) \right) \right| = 1 \; ,
\]
this change of variables is volume preserving. Thus, this implies that
\[
\omega( U_t \varphi, U_t \Phi ) = \int_{\mathbb R^{2n_x}} \nabla \varphi( z' )^\tr \Omega \nabla \Phi( z' ) \, \mathrm{d}z' = \omega(\varphi,\Phi)
\]
for all $\varphi,\Phi \in \mathbb W_{1,2}^{2n_x}$. Thus, $U_t$ is a symplectic operator.
\qed

\subsection{Symplectic duality}
\label{sec::duality}
In order to prepare the analysis of the properties of the spectrum of symplectic Koopman operators, it is helpful to introduce a notion of duality. However, instead of defining duality in a Hilbert space, we propose to introduce the following notion of symplectic duality.
\begin{defn}
Let $A: \mathbb W_{1,2}^{2n_x} \to \mathbb W_{1,2}^{2n_x}$ be a linear operator. If there exists a linear operator \mbox{$A^\star: \mathbb W_{1,2}^{2n_x} \to \mathbb W_{1,2}^{2n_x}$} with
\[
\omega( \varphi, A \Phi ) = \omega( A^\star \varphi, \Phi )
\]
for all $\varphi,\Phi \in \mathbb W_{1,2}^{2n_x}$, then we call $A^\star$ the symplectic adjoint operator of $A$.
\end{defn}

\noindent
In the following, we consider a linear differential operator $L: \mathbb W_{1,2}^{2n_x} \to \mathbb W_{1,2}^{2n_x}$ that is given by
\begin{align}
\label{eq::L}
L \Phi = F^\tr \nabla \Phi
\end{align}
for all $\Phi \in \mathbb W_{1,2}^{2n_x}$. Notice that this operator is found by differentiating the definition of the Pontryagin-Koopman operator $U_t$ with respect to time, such that $L = \dot U_t$.

\begin{lem}
Let Assumptions~\ref{ass::f},~\ref{ass::l}, and~\ref{ass::u} hold. The operator $L$ admits a symplectic adjoint $L^\star$ in the symplectic space $(\mathbb W_{2,2}^{2n_x}, \omega)$. This adjoint is given by
\begin{align}
\label{eq::LL}
L^\star = -L \; .
\end{align}
\end{lem}

\textbf{Proof.}
For any function $\Phi \in \mathbb W_{2,2}^{2n_x}$, we have
\begin{align}
\label{eq::NablaL}
\nabla (L \Phi) = \frac{\mathrm{d}F}{\mathrm{d}y} \nabla \Phi + \nabla^2 \Phi F \; ,
\end{align}
which follows by differentiating~\eqref{eq::L} on both sides. Next, we know from Theorem~\ref{thm::symplectic} that $U_t$ is a symplectic operator. Thus, we can differentiate the equation
\[
\omega( U_t \varphi, U_t \Phi ) = \omega( \varphi, \Phi )
\]
on both sides with respect to $t$, which yields
\[
\frac{\mathrm{d}}{\mathrm{d}t}\omega( U_t \varphi, U_t \Phi ) = 0 \; .
\]
Now, if $\varphi,\Phi \in \mathbb W_{2,2}^{2n_x}$, we can expand the derivative on the right as
\begin{align}
0 = \frac{\mathrm{d}}{\mathrm{d}t} \, \omega( U_t \varphi, U_t \Phi ) \notag 
&\overset{\eqref{eq::aux4}}{=} \frac{\mathrm{d}}{\mathrm{d}t} \, \int_{\mathbb R^{2n_x}} \nabla \varphi( \Gamma_t(z) )^\tr \Omega \nabla \Phi( \Gamma_t(z) ) \, \mathrm{d}z \notag \\
&= \int_{\mathbb R^{2n_x}}  F(\Gamma_t(z))^\tr \nabla^2 \varphi( \Gamma_t(z) ) \Omega \nabla \Phi( \Gamma_t(z) ) \, \mathrm{d}z 
 + \int_{\mathbb R^{2n_x}} \nabla \varphi( \Gamma_t(z) )^\tr \Omega \nabla^2 \Phi( \Gamma_t(z) ) F(\Gamma_t(z)) \, \mathrm{d}z \notag
\end{align}
and, after changing variables, $z \leftarrow \Gamma_t(z)$, recalling that $\mathrm{det}\left( \frac{\partial}{\partial z} \Gamma_t \right) = 1$, and resorting terms,
\begin{equation}
\int_{\mathbb R^{2n_x}} F(z)^\tr \nabla^2 \varphi(z) \Omega \nabla \Phi(z) \, \mathrm{d}z 
\label{eq::aux5}
= -\int_{\mathbb R^{2n_x}} \nabla \varphi(z)^\tr \Omega \nabla^2 \Phi(z) F(z) \, \mathrm{d}z \; .
\end{equation}
Moreover, we recall that the equation
\begin{align}
\label{eq::aux6}
\Omega \frac{\mathrm{d}F}{\mathrm{d}y} = \left( \frac{\mathrm{d}F}{\mathrm{d}y}\right)^\tr \Omega^\tr = - \left( \frac{\mathrm{d}F}{\mathrm{d}y}\right)^\tr \Omega
\end{align}
also holds (see Equation~\eqref{eq::FHamiltonian}). Thus, in summary, we have
\begin{align*}
\omega(\varphi , L \Phi ) = \displaystyle\int_{\mathbb R^{2n_x}} \left[  \nabla \varphi^\tr \Omega \nabla [ L \Phi ]  \right](z) \, \mathrm{d}z  
&\quad \overset{\eqref{eq::NablaL}}{=} \qquad \displaystyle\int_{\mathbb R^{2n_x}} \left[  \nabla \varphi^\tr \Omega \left(  \frac{\mathrm{d}F}{\mathrm{d}x} \nabla \Phi + \nabla^2 \Phi F \right)  \right](z) \, \mathrm{d}z \notag \\[0.2cm]
&\overset{\eqref{eq::aux5},\eqref{eq::aux6}}{=} - \displaystyle\int_{\mathbb R^{2n_x}} \left[  \left(  \frac{\mathrm{d}F}{\mathrm{d}x} \nabla \varphi + \nabla^2 \varphi F \right)^\tr \Omega \nabla \Phi   \right](z) \, \mathrm{d}z \notag \\[0.2cm]
&\quad \, = \qquad \omega(-L \varphi , \Phi ) \; .
\end{align*}
Because the latter equation holds for all $\varphi,\Phi \in \mathbb W_{2,2}^{2n_x}$, we find that $L^\star = - L$ is the symplectic adjoint of $L$, as claimed by the statement of this lemma.
\qed

\begin{rem}
In many articles on Koopman operators, a so-called transport differential equation of the form
\[
\frac{\partial}{\partial t} \phi = - \mathrm{div}(F \phi)
\]
is considered. If one introduces suitable boundary conditions, this advection PDE can be interpreted as the derivative of a Perron-Frobenius operator that is dual to the Koopman operator $U_t$ with respect to the standard $\mathbb L_2$-scalar product~\cite{Ding1998}. Notice that in our case, the symplectic adjoint operator
\[
L^\star: \; \phi \, \to \, - F^\tr \nabla \phi = -\mathrm{div}(F \phi) + \underbrace{\mathrm{div}(F)}_{=0} \phi
\]
happens to coincide with the right-hand operator of the above advection PDE. Thus, physically, one could interpret the operator $L^\star$ as an advection operator of the incompressible flow field $F$ recalling that $\mathrm{div}(F)=0$.
\end{rem}

\section{Spectral Analysis}
\label{sec::spectrum}
In this section, we show that the symplectic structures of the Koopman operator, as analyzed in the previous section, have important consequences on its set of eigenvalues. Moreover, Section~\ref{sec::Manifolds} uses these spectral properties of the Koopman operator to characterize global optimal control laws that are associated with the infinite horizon optimal control problem~\eqref{eq::ocp}.

\subsection{Eigenfunctions and Eigenvalues}
As already mentioned in the introduction, the spectrum of general Koopman operators has been analyzed by many authors; for example in~\cite{Arbabi2017,Budisic2012,Mezic2005}. In the context of this paper, we call a weakly differentiable function $\Psi$ an eigenfunction of the differential Pontryagin-Koopman operator $L$, if the Lebesgue measure of the set $\, \mathrm{supp}(\Psi) \,$ is either unbounded, or, if it exists, is not equal to $0$, and
\[
L \Psi = \kappa \Psi \qquad \Longleftrightarrow \qquad \nabla \Psi^\tr F = \kappa \Psi
\]
for a potentially complex eigenvalue $\kappa \in \mathbb C$. The fact that the symplectic adjoint of the operator $L$ is given by $L^\star = -L$ has important consequences on its spectrum, which can be summarized as follows.

\begin{thm}
\label{thm::mirrorEigenvalues}
Let Assumptions~\ref{ass::f},~\ref{ass::l}, and~\ref{ass::u} hold and let $\Psi_1,\Psi_2 \in \mathbb W_{2,2}^{2n_x}$ be two eigenfunctions of $L$ with eigenvalues $\kappa_1,\kappa_2 \in \mathbb C$. If $\omega(\Psi_1,\Psi_2) \neq 0$, then we must have $\kappa_1 = -\kappa_2$.
\end{thm}

\textbf{Proof.}
Because the functions $\Psi_1,\Psi_2 \in \mathbb W_{2,2}^{2n_x}$ are eigenfunctions of $L$, they satisfy
\begin{align}
\label{eq::aux7}
L \Psi_1 = \kappa_1 \Psi_1 \quad \text{and} \quad L \Psi_2 = \kappa_2 \Psi_2 \; .
\end{align}
Thus, since $\omega$ is a bilinear form, we have
\begin{equation*}
\kappa_2 \omega(\Psi_1, \Psi_2) \overset{\eqref{eq::aux7}}{=} \omega( \Psi_1, L \Psi_2) 
\overset{\eqref{eq::LL}}{=} -\omega( L \Psi_1, \Psi_2)
\overset{\eqref{eq::aux7}}{=} -\kappa_1 \omega(\Psi_1, \Psi_2) \notag
\end{equation*}
and, after resorting terms,
\[
(\kappa_1 + \kappa_2 )\omega(\Psi_1, \Psi_2) = 0 \; .
\]
Thus, if $\omega(\Psi_1, \Psi_2) \neq 0$, we must have $\kappa_1 + \kappa_2 = 0$. This is equivalent to the statement of the theorem.
\qed

\begin{rem}
The statement of Theorem~\ref{thm::mirrorEigenvalues} is formally not directly applicable for eigenfunctions $\Psi_1,\Psi_2$ of $L$ that are locally weakly twice differentiable but whose derivatives are not square integrable, such that $\Psi_1,\Psi_2$ are not elements of the Sobolev space~$\mathbb W_{2,2}^{2n_x}$. However, in practice, one is usually interested in constructing eigenfunctions of $L$ on a compact domain $C \subset \mathbb R^{2n_x}$, in the following called the region of interest, such that
\begin{align}
\label{eq::compact}
\forall x \in C, \quad [L \Psi](x) = \kappa \Psi(x) \; .
\end{align}
This region of interest $C$ can, for example, model a-priori bounds on the optimal primal and dual trajectories of~\eqref{eq::ocp}. Consequently, because we are simply not interested in how the flow $\Gamma_t$ and the associated eigenfunctions $\Psi$ are defined outside of the set $C$, these functions can simply be redefined arbitrarily for $x \notin C$, for example, such that the desired eigenfunctions of $L$ satisfy $\Psi \in \mathbb W_{2,2}^{2 n_x}$ by construction. Thus, for the purpose of this paper, it is not restrictive at all to assume that the derivatives of the eigenfunctions of $L$ are square integrable. Notice that this technique is also illustrated by our tutorial case study in Section~\ref{sec::numerics}, where we explain how to choose $C$ and how to discretize~\eqref{eq::compact} on $C$.
\end{rem}

Let $\sigma(L) \subseteq \mathbb C$ denote the spectrum of $L$; that is, the set of eigenvalues of the linear operator $L$. In the following, we use the notation $\Psi_{\kappa} \in \mathbb W_{2,2}^{2n_x}$ to denote an eigenfunction that is associated with an eigenvalue $\kappa \in \sigma(L)$. Moreover, for a function $q: \mathbb R^{2 n_x} \to  \mathbb R^{2 n_x}$, we use the shorthand notation
\[
\omega(q,q) = \left(
\begin{array}{cccc}
\omega(q_1,q_1) & \omega(q_1,q_2) & \ldots & \omega(q_1,q_{2n_x}) \\
\omega(q_2,q_1) & \omega(q_2,q_2) & \ldots & \omega(q_2,q_{2n_x}) \\
\vdots & \vdots & \ddots & \vdots \\
\omega(q_{2n_x},q_1) & \omega(q_{2n_x},q_2) & \ldots & \omega(q_{2n_x},q_{2n_x}) \\
\end{array}
\right)
\]
to denote the matrix that is obtained by evaluating the skew symmetric bilinear form $\omega$ for all possible combinations of the components of $q$. We call $\, q \,$ a skew-orthogonal function in the symplectic space $(\mathbb W_{1,2}^{2n_x}, \omega)$, if it satisfies $\omega(q,q) = I$. Notice that the construction of such skew orthogonal functions is straightforward by using the standard skew-symmetric variant of the Gram-Schmidt algorithm~\cite{Arnold2001}.

\begin{defn}
\label{def::spectralUnit}
The operator $L$ is said to admit a spectral decomposition with respect to a given function $q$, if there exists a generalized function \mbox{$a: \sigma(L) \to \overline{\mathbb C}^{2 n_x}$} such that
\begin{align}
\label{eq::spectralDecomposition}
q = \int_{\sigma(L)} \Psi_{\kappa} a(\kappa) \, \mathrm{d}\kappa \; .
\end{align}
\end{defn}

Notice that the above definition uses the ``control engineering notation'' for generalized functions, which means that we use $a$ as if it was a standard function, although this notation suppresses the distributional nature of $a$. Thus, in mathematical terms, this notation has to be translated as ``$a$ represents a distribution of order $0$; that is, a linear operator on $\sigma(L)$ that is Lipschitz continuous with respect to the $L_\infty$-norm''.\footnote{Notice that, for example, $a$ could be a Dirac distribution, which is not a function in the traditional sense but, by construction, a Lipschitz continuous linear operator.}
The following statement is a consequence of Theorem~\ref{thm::mirrorEigenvalues}.

\begin{cor}
\label{cor::Kspectrum}
Let Assumptions~\ref{ass::f},~\ref{ass::l}, and~\ref{ass::u} hold.
Let the operator $L$ admit a spectral decomposition with respect to a skew-orthogonal
function $q$. Then, there exist at least $2n_x$ eigenvalues
\[
\kappa_1^+, \kappa_2^+, \ldots, \kappa_{n_x}^+, \kappa_1^-, \kappa_2^-, \ldots, \kappa_{n_x}^- \in \sigma(L) \; ,
\]
such that $\kappa_i^+ = -\kappa_i^{-}$, where $\kappa_i^+$ has a non-negative real part for all $i \in \{ 1, 2, \ldots, n_x \}$.
\end{cor}

\textbf{Proof.}
By substituting~\eqref{eq::spectralDecomposition} in the equation $\omega(q,q) = I$, we find that the equation
\begin{equation}\label{eq::auxspect}
I = \int_{\sigma(L)} \int_{\sigma(L)} \omega( \Psi_{\kappa'}, \Psi_{\kappa} ) \underbrace{ a(\kappa') a(\kappa)^\tr }_{\mathrm{rank} \; 1} \, \mathrm{d}\kappa \, \mathrm{d}\kappa' \; ,
\end{equation}
holds. Let us have a closer look at the terms in~\eqref{eq::auxspect}. Clearly, the unit matrix on the left has full rank. But, on the other side, we have an integral over the rank $1$ matrices $a(\kappa') a(\kappa)^\tr$. This integral term can only have full rank, if there are at least $2n_x$ pairs of eigenvalues $(\kappa',\kappa) \in \sigma(L) \times \sigma(L)$ for which
$$\omega(\Psi_{\kappa'}, \Psi_{\kappa}) \neq 0 \; .$$
But now Theorem~\ref{thm::mirrorEigenvalues} implies that all these pairs must be such that $\kappa' = -\kappa$ and, after sorting all eigenvalues with respect to their real-part, we find that there must be at least $n_x$ eigenvalues with non-negative real part and $n_x$ associated mirrored eigenvalues with non-positive real-part, as claimed by the statement of this corollary.
\qed

\begin{rem}
Notice that Corollary~\ref{cor::Kspectrum} makes a statement about the spectrum of $L$ under the assumption that this operator admits a spectral decomposition with respect to at least one skew orthogonal function. Thus, one should further ask the question under which assumptions the existence of such a spectral decomposition can be guaranteed for at least one such skew orthogonal function. In full generality, this question is difficult to answer, but sufficient conditions for the existence of (much more general and, in our context, sufficient) spectral decompositions can be found in~\cite{Arbabi2017,Budisic2012,Mauroy2016,Mezic2005}, which use ideas from the field of ergodic theory~\cite{Wiener1941} as well as Yoshida's theorem~\cite{Yosida1978}. For example, if the monodromy matrix, $\frac{\partial}{\partial z} \Gamma_T( x_\mathrm{p}(0) )$, of a periodic optimal orbit with period length $T > 0$ is diagonalizable, one can ensure that a spectral decomposition is possible. This sufficient condition follows by applying Proposition~3 in~\cite{Mauroy2016} to the Pontryagin differential equation. A more complete review of such results from the field of functional analysis would, however, go beyond the scope of the present paper.
\end{rem}

\section{Optimal Feedback Control Laws}
\label{sec::Manifolds}

The theoretical results from the previous sections can be used to derive optimality conditions, which, in turn, can be used to develop practical numerical algorithms for solving~\eqref{eq::ocp}. Let us introduce the set
\begin{align}
\sigma^+(L) &= \{ \kappa \in \sigma(L) \mid \mathrm{Re}(\kappa) > 0 \} \notag
\end{align}
of unstable eigenvalues and its associated invariant manifold
\begin{align}
\mathcal M^+ &= \left\{ [x^\tr,\lambda^\tr]^\tr \in \mathbb R^{2n_x} \; \middle| \; \begin{array}{l}
\forall \kappa \in \sigma^+(L), \\
\Psi_{\kappa}(x,\lambda) = 0                                                                            \end{array}
\right\} \; . \notag
\end{align}
Because we assume that the eigenfunctions $\Psi_{\kappa}$ are weakly differentiable, we may assume without loss of generality that the functions $\Psi_{\kappa}$ are also \mbox{continuous---otherwise} there exists a continuous function $\tilde \Psi_{\kappa}(x) = \Psi_{\kappa}(x)$ for almost every $x \in \mathbb R^{2n_x}$ and we can use $\tilde \Psi_{\kappa}$ instead of $\Psi_{\kappa}$ recalling that we work with Sobolev spaces, in which such arguments are indeed possible.

In the following, we say that $\mu^\star$ is a regular optimal control law of~\eqref{eq::ocp}, if the closed-loop trajectories,
\[
\dot x^\star(t) = f( x^\star(t), \mu^\star( x^\star(t) ) ) \quad \text{with} \quad x^\star(0) = x_0 \; ,
\]
are minimizers of~\eqref{eq::ocp} at which Pontryagin's necessary optimality conditions are satisfied such that $x^\star$ converges to an optimal periodic orbit at which the conditions of Proposition~\ref{prop::multiplier} are satisfied.

\begin{thm}
\label{thm::OptimalControl}
Let Assumptions~\ref{ass::f},~\ref{ass::l}, and~\ref{ass::u} hold. We further assume (without loss of generality) that the eigenfunctions $\Psi_{\kappa} \in W_{2,2}^{2 n_x}$ are continuous. Let \mbox{$\mu^\star: \mathbb R^{n_x} \to \mathbb R^{n_x}$} be a regular optimal control law of~\eqref{eq::ocp}. Then, there exists a function \mbox{$\Lambda: \mathbb R^{n_x} \to \mathbb R^{n_x}$} such that
\[
[ x^\tr , \Lambda(x)^\tr ]^\tr \in \mathcal M^+ \quad \text{and} \quad \mu^\star(x) = u^\star( x, \Lambda(x)) \; .
\]
\end{thm}

\textbf{Proof.}
Since $L$ is a time-autonomous infinitesimal generator of $U_t$, the eigenfunctions of $L$ also satisfy~\cite{Mezic2005}
\begin{align}
\label{eq::Ut}
U_t \Psi_{\kappa} = e^{\kappa t} \Psi_{\kappa}
\end{align}
for all $\kappa \in \sigma(L)$. The remainder of the proof is divided into two parts. In the first part, we show that~\eqref{eq::Ut} implies that the optimal periodic limit orbit is a subset of the manifold $\mathcal M^+$ and, in the second part, we use this property of the limit orbit to construct a globally optimal feedback law.

\smallskip
\textit{Part I:}
Let $(x_{\mathrm{p}},\lambda_{\mathrm{p}})$ denote the optimal periodic limit orbit of~\eqref{eq::ocp} with period time $T$. If we would have \mbox{$(x_{\mathrm{p}}(t),\lambda_{\mathrm{p}}(t)) \notin \mathcal M^+$} for a time $t \in [0,T]$, we must have
\[
\Psi_{\kappa}(x_{\mathrm{p}}(t),\lambda_{\mathrm{p}}(t)) \neq 0
\]
for at least one $\kappa \in \sigma^+(L)$. Because the optimal periodic orbit satisfies Pontryagin's differential equation, this implies in particular that
\begin{equation}
\left| \Psi_{\kappa}(x_{\mathrm{p}}(t+T),\lambda_{\mathrm{p}}(t+T)) \right| \overset{\eqref{eq::Ut}}{=} \left| e^{\kappa T} \Psi_{\kappa}(x_{\mathrm{p}}(t),\lambda_{\mathrm{p}}(t)) \right| \notag 
> \left| \Psi_{\kappa}(x_{\mathrm{p}}(t),\lambda_{\mathrm{p}}(t)) \right| \,
\end{equation}
as $\kappa$ has a strictly positive real part. But this is a contradiction, since $(x_{\mathrm{p}},\lambda_{\mathrm{p}})$ is periodic. Thus, in summary, we must have 
\begin{align}
\label{eq::mp}
\forall t \in [0,T], \quad (x_{\mathrm{p}}(t),\lambda_{\mathrm{p}}(t)) \in \mathcal M^+ \; .
\end{align}

\smallskip
\textit{Part II:} Let $(x,\lambda)$ denote an optimal solution of~\eqref{eq::ocp} with $x(0) = x_0$. Now, if we would have $(x(0),\lambda(0)) \notin \mathcal M^+$, then there would have to exist at least one $\kappa \in \sigma^+(L)$ for which
\[
\Psi_{\kappa}( x(0), \lambda(0) ) \neq 0 \; .
\]
But if this would be the case, then we would also have
\begin{align}
\lim_{t \to \infty} |\Psi_{\kappa}( x(t), \lambda(t) )| = \infty \; ,
\end{align}
as $\Psi_{\kappa}$ is strictly unstable. But this limit statement is in conflict with~\eqref{eq::mp}, since $\Psi_{\kappa}$ is continuous and we assume that $(x(t),\lambda(t))$ converges to the optimal periodic limit orbit. Thus, in summary, there exists for every $x_0$ a $\lambda_0 = \lambda(0)$ with $(x_0,\lambda_0) \in \mathcal M^+$ and the corresponding map from $x_0$ to $\lambda_0$ can be denoted by $\Lambda$. The associated optimal control input is given by
\[
\mu^\star(x_0) = u^\star(x_0,\Lambda(x_0)) \; .
\]
This is already sufficient to establish the statement of this theorem, as the optimal feedback law must be time-autonomous.\qed

Notice that Theorem~\ref{thm::OptimalControl} can be used to systematically search for globally optimal solutions of~\eqref{eq::ocp}. Here, one of the key observations is that if, in addition to the assumptions of Theorem~\ref{thm::OptimalControl}, the assumptions of Corollary~\ref{cor::Kspectrum} are also satisfied and if none of the non-negative eigenvalues $\kappa_1^+, \kappa_2^+, \ldots, \kappa_{n_x}^+$ happens to be on the imaginary axis, then the parametric equation system
\begin{align}
\label{eq::Lambda}
\forall \kappa \in \sigma^+(L), \quad \Psi_{\kappa}(x_0,\lambda_0 ) = 0 \; .
\end{align}
consists of at least $n_x$ independent equations while the number of variables, $\mathrm{dim}(\lambda) = n_x$, is equal to $n_x$, too. Thus, this equation system can, in many practical instances, be expected to admit a finite number of parametric solutions \mbox{$\lambda_0 = \Lambda(x_0)$} only. This means that, if all the above regularity assumptions are satisfied, Theorem~\ref{thm::OptimalControl} singles out a finite number of candidate control laws \mbox{$\mu^\star(x) = u^\star(x,\Lambda(x))$---at} least one of which must be globally optimal.

\section{Numerical example}
\label{sec::numerics}

This section illustrates how Theorem~\ref{thm::OptimalControl} can be used to construct accurate approximations
of globally optimal control laws of~\eqref{eq::ocp} for a Van der Pol oscillator system with
\begin{align}
\label{eq::oscillator}
f(x,u) &= \left(
\begin{array}{l}
x_2 \\
x_1 - \frac{1}{2}\left( 1-x_1^2 \right) x_2 + x_1 u
\end{array}
\right) \quad \\
\text{and} \quad l(x,u) &=  \frac{1}{2}\left( x_2^2 + u^2 \right) \; .
\end{align}
The associated system, $\dot x(t) = f(x(t),u(t))$, has a linearly uncontrollable
equilibrium at $x = (0,0)^\intercal$. Notice that for this particular example an explicit
solution of the Hamilton-Jacobi-Bellman equation is known~\cite{Rodrigues2016}:
it turns out that the optimal value function $V$ and the globally optimal feedback law $\mu^\star$
are given, respectively, by
\begin{align}
\label{ex::exactSolution}
V(x) = \frac{1}{2}( x_1^2 + x_2^2 ) \quad \text{and} \quad \mu^\star(x) = -x_1 x_2 \; .
\end{align}
In the following, this explicit solution is, however, only used to assess the accuracy of the
proposed method; that is, the numerical procedure below neither knows the above expression
for $V$ nor for $\mu^\star$.

\subsection{Galerkin discretization}
In order to illustrate the practical applicability of the developments in Section~\ref{sec::Manifolds}, we introduce
a simple Galerkin discretization of the operator $L$. Let $\varphi_{1},\ldots,\varphi_{N}\in 
\mathbb W_{1,2}^{2n_x}$
be orthogonal functions with respect to the standard $\mathbb L_2$-scalar product $\langle \cdot, \cdot \rangle$ on $\mathbb W_{1,2}^{2n_x}$. Next, if we compute the coefficients
\begin{gather*}
M_{i,j} = \langle L \varphi_i, \varphi_j \rangle \;  
\end{gather*}
for all $i,j \in\{1,\ldots,N\}$, the matrix $M^\tr$ can be interpreted as a Galerkin
approximation of the operator $L$ over the subspace spanned
by $\varphi_{1},\ldots,\varphi_{N}$. In our implementation, we set $\varphi_i$ to the $i$th multivariate Legendre polynomial on the $4$-dimensional compact interval box $C = [-\frac{1}{2},\frac{1}{2}]^4$ and we set $\varphi_i(x,\lambda) = 0$ outside of this domain; that is, for $(x,\lambda) \notin C$. Consequently, our discretization can only be expected to be accurate inside our region of interest $C$, but other choices for $C$ and for the basis functions would be possible, too.

\begin{rem}
Standard Galerkin methods are, in general, numerically unstable when applied to advection operators and, consequently, although this method happens to yield reasonable approximations for our particular example, such naive discretization schemes cannot be recommended in general. More advanced discretization schemes for linear advection operators can be found in the modern PDE literature; see, for example~\cite{Roos2008}. A more complete discussion of such numerical discretization methods for Pontryagin-Koopman operators, is, however, beyond the scope of this paper.
\end{rem}

\begin{figure*}[t]
\centering
\includegraphics[width=0.33\textwidth]{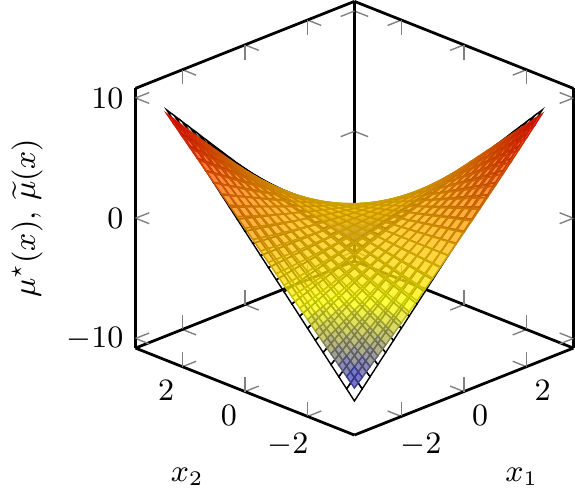}%
\hspace{0.5em}%
\includegraphics[width=0.3\textwidth]{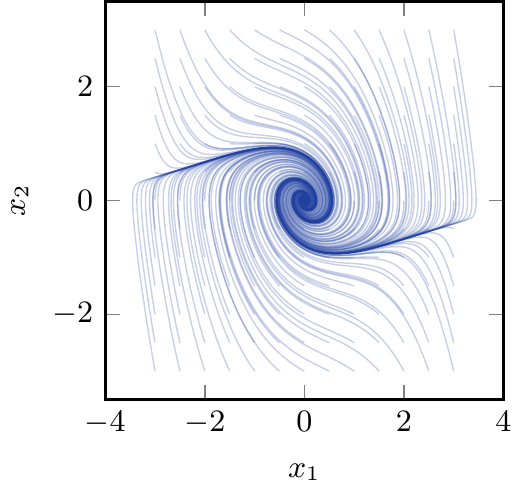}%
\hspace{0.5em}%
\includegraphics[width=0.3\textwidth]{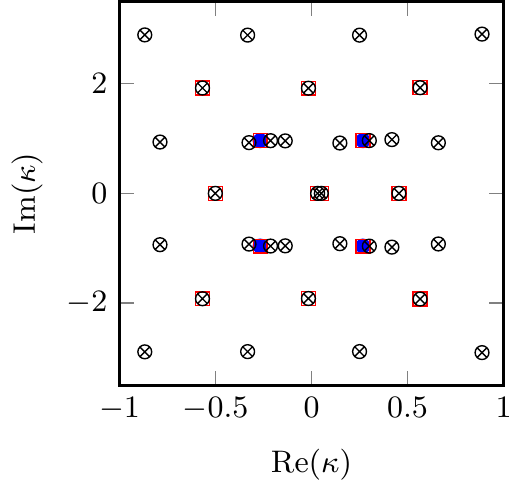}
\caption{\label{fig::closedloop} Left: Comparison between the optimal feedback law 
$\mu^\star(x) = -x_1x_2$, and the approximate feedback law 
$\widetilde{\mu}(x)$ in~\eqref{eq::mutilde}. Middle: trajectories of the 
closed-loop system under the approximate feedback law. Right: Eigenvalues of 
the projection $M \in \mathbb R^{N \times N}$ for $N=4$ (blue circles), $N=15$ (red squares),
and $N=35$ (black circles with a cross).} 
 \end{figure*}

\subsection{Approximations of the optimal feedback law}
The above Galerkin approximation of the operator $L$ can be used to construct approximate eigenfunctions. In detail, if 
$a\in\mathbb{R}^{N}$ is a left eigenvector of $M$ with eigenvalue $\kappa \in \mathbb C$, then
$$\Psi = \sum_{i=1}^N a_i \varphi_i \quad \Longrightarrow \quad L\Psi \approx \kappa \Psi$$
is an approximation of an eigenfunction of $L$. The right plot in Figure~\ref{fig::closedloop} shows the spectrum of the matrix $M$ for different choices of $N$ (blue circles: $N=4$, red squares: $N=15$, and black circles with a cross: $N=35$). In order to understand the structure of this spectrum it is helpful to recall  Corollary~\ref{cor::Kspectrum}, which predicts that there exits at least $2$ eigenvalues $\kappa_1,\kappa_2$ of $L$ such that $-\kappa_1$ and $-\kappa_2$ are also eigenvalues of $L$. Notice that such symmetric eigenvalue pairs are indeed present in the spectrum of $M$, although $M$ is only a Galerkin approximation of $L$.

In order to further illustrate how the above spectral analysis of $M$ can be used to construct approximations of the globally optimal control law $\mu^\star$, we can compute an approximation of the manifold $\mathcal M^+$ by using the approximate eigenfunctions instead of the exact ones (see Theorem~\ref{thm::OptimalControl}). For example, for $N=4$, the matrix $M$ has the eigenvalues
\begin{equation}
\kappa^{\pm}_{1,2} = \frac{1}{48}(\beta_{m} \pm \beta_{p} \sqrt{-1}) \quad
\text{with} \quad \beta_{p} = \sqrt{ +983 + 96 \sqrt{143} } \quad
\text{and} \quad \beta_{m} = \sqrt{ -983 + 96 \sqrt{143} }
\end{equation}
and the associated Galerkin approximation of the globally optimal control law is given by
\begin{equation}
\label{eq::mutilde}
\widetilde{\mu}(x) = 
(12 - \sqrt{143}) x_1^2 + \frac{(11 - \beta_m )}{2} x_1x_2 \;.
\end{equation}
The left plot in Figure~\ref{fig::closedloop} shows $\widetilde \mu$ and compares it to the optimal feedback law $\mu^\star$. In fact, the squared integral error over $C$ is 
approximately $6\times 10^{-5}$. Quite remarkably, this approximate feedback law can even be used to control the system for initial values outside of $C$. The plot in the middle of Figure~\ref{fig::closedloop} shows the corresponding trajectories of the closed-loop system that are obtained by using the approximately optimal feedback law $\widetilde \mu$.

\section{Conclusions}
\label{sec::conclusion}

This paper has presented an analysis of infinite horizon nonlinear optimal control
problems, whose minimizers satisfy Pontryagin's necessary conditions of optimality.
The proposed formalism is based on Pontryagin-Koopman operators, which have
been shown to possess a symplectic structure, as revealed by Theorem~\ref{thm::symplectic}.
Moreover, Theorem~\ref{thm::mirrorEigenvalues} and Corollary~\ref{cor::Kspectrum} have
established conditions under which the spectrum of the differential Pontryagin-Koopman
operator contains at least $2n_x$ mirrored eigenvalues. This spectral structure is used in
Theorem~\ref{thm::OptimalControl} to characterize optimal control laws.

The theoretical findings of this paper have been applied to construct accurate
approximations of a globally optimal control law for a Van der Pol oscillator, which
illustrates the potential of the proposed Pontryagin-Koopman operator based framework
for the design of global optimal control algorithms. Here, it needs to be highlighted that,
in contrast to dynamic programming methods, which rely on the discretization of \textit{nonlinear}
Hamilton-Jacobi-Bellman PDEs, the proposed framework for global
optimal control relies on the computation of eigenfunctions of a \textit{linear} differential operator.
This opens the door to an application of linear algebra methods and tailored discretization schemes
for linear PDEs, which have never been considered for the computation of optimal control laws.
Therefore, the development of tailored, structure-exploiting projection and linear algebra
methods for symplectic Pontryagin-Koopman operators and their application to
global optimal control can be regarded as a promising direction for future research.

\bibliographystyle{plain}
\bibliography{references}

\end{document}